\documentclass[12pt,twoside]{article}
\usepackage{amsmath}
\usepackage{mathrsfs}
\usepackage{amsfonts}
\usepackage{amssymb}

\usepackage{amsfonts,amsmath,amsthm}
\allowdisplaybreaks

\numberwithin{equation}{section} \textwidth 15.5 true cm
\begin{document}
  \title{{\bf Global existence and optimal decay rates for three-dimensional compressible viscoelastic flows }}
  \author{Xianpeng Hu$^a$, Guochun Wu$^b$ \thanks{Corresponding author: Guochun Wu, Email
 address: guochunwu@126.com}\\{\it\scriptsize $^a$Courant Institute of Mathematical Sciences, New York University,
                           New York, NY 10012, USA.}\\{\it\scriptsize $^b$School of
Mathematical Sciences, Xiamen University,
 Fujian 361005, China}\\{\it\small
 Email
 address: xianpeng@cims.nyu.edu, guochunwu@126.com}}
  \date{}
  \maketitle
  \begin{abstract}{\small In this paper, we are concerned with the global existence and
  optimal rates of strong solutions for three-dimensional compressible
viscoelastic flows. We prove the global existence of the strong
solutions by the standard energy method under the condition that
the initial data are close to the constant equilibrium state in
$H^2$-framework. If additionally the initial data belong to $L^1$,
the optimal convergence rates of the solutions in $L^p$-norm with
$2\leq p\leq 6$ and optimal convergence rates of their spatial
derivatives in $L^2$-norm are obtained.
 }
  \end{abstract}
\section{Introduction}

$\ \ \ \ $In this paper, we are interested in three-dimensional
compressible viscoelastic flows [3,5,12,23]:
$$
\rho_t+div (\rho u)=0,\eqno(1.1a)$$$$ (\rho u)_t+div(\rho u\otimes
u)-\mu \triangle u-(\lambda+\mu)\nabla div u+\nabla P(\rho)=\alpha
div (\rho FF^{T}),\eqno(1.1b)$$$$ F_t+u\cdot\nabla F=\nabla u
F,\eqno(1.1c)$$ for $(t,x)\in[0,+\infty)\times\mathbb{R}^{3}$.
Here $\rho$, $u\in\mathbb R^3$, $F\in M^{3\times 3}$ (the set of
$3\times 3$ matrices with positive determinants) denote the
density, the velocity, and the deformation gradient, respectively.
The Lam\'e coefficients $\mu$ and $\lambda$ are satisfied the
physical condition:
$$\mu>0,\ \ 2\mu+3\lambda> 0,$$
which ensures that the operator $-\mu\Delta-(\lambda+\mu)\nabla div$
is a strongly elliptic operator. The pressure term $P(\rho)$ is an
increasing and convex function of $\rho$ for $\rho>0$. The symbol
$\otimes$ denotes the Kronecker tensor product, $F^T$ means the
transpose matrix of $F$, and the notation $u\cdot\nabla F$ is
understood to be $(u\cdot\nabla)F$. For system (1.1), the
corresponding elastic energy is chosen to be the special form of the
Hookean linear elasticity:
$$W(F)=\frac{\alpha}{2}|F|^2+\frac{1}{\rho}\int_0^{\rho}P(s)ds,\ \ \alpha>0,$$
which, however, does not reduce the essential difficulties for
analysis. Indeed, all the results we describe here can be
generalized to a more general cases.

In this paper, we investigate the Cauchy problem of system (1.1)
with the initial condition:
$$(\rho,u,F)|_{t=0}=(\rho_0(x),u_0(x),F_0(x)),\ \ x\in\mathbb R^3.\eqno(1.2)$$
We also assume that $$div(\rho F^T)=0,\
F^{lk}(0)\nabla_lF^{ij}(0)=F^{lj}(0)\nabla_lF^{ik}(0).\eqno(1.3)$$
It is standard that the condition (1.3) is preserved by the flow,
which has been proved in [7,21].

For the incompressible viscoelastic flows and related models,
there are many important progress on classical solutions, refer to
[1,2,9,13,16] and references therein. On the other hand, the
global existence of weak solutions to the incompressible
viscoelastic flows with large initial data is still an outstanding
open question, although there are some progress in that direction
[15,17,18]. For the compressible viscoelastic flows, to our
knowledge, there are few results on the dynamics of global
solutions to compressible viscoelastic flows, especially on the
large time behavior. The local existence of multi-dimensional
strong solution was obtained in [6], and the global existence of
strong solution with the lowest regularity was shown in [7,21].
For the initial boundary value problem, global in time solution
was proved to exist uniquely near the equilibrium state in [8,22].

In this paper, we firstly study the optimal time-decay rate of the
global strong solutions to the Cauchy problem (1.1)-(1.2). To be
more precise, the main purpose of this paper is to study the
existence and uniqueness of global strong solutions and in
particular the asymptotic behavior on the Cauchy problem of
compressible viscoelastic flows. We prove the global existence of
strong solutions by the standard energy method in spirit of
Matsumura and Nishida [19,20]. In order to obtain the linear
time-decay estimates, we need to analysis the properties of the
semigroup, as in [10,11,14,24]. Unfortunately, it seems untractable,
since the system (1.1) has thirteen equations. To overcome this
difficulty, we take Hodge decomposition of the linear system, then
it becomes two similar systems, each of those only involves two
variables, which makes us be able to obtain the optimal time-decay
estimates.

Our main results are formulated in the following theorem:\\
\\{\bf Theorem 1.1.} Assume that the initial value $(\rho_0-1,u_0,F_0-I)\in H^2(\mathbb R^3)$
satisfies the constraints (1.3), then there exists a constant
$\delta_0$ such that if $$|(\rho_0-1,u_0,F_0-I)|_{H^2}\leq
\delta_0,\eqno(1.4)$$ then there exists a unique globally strong
solution $(\rho,u,F)$ of the Cauchy problem $(1.1)-(1.2)$ such
that for any $t\in [0,\infty)$,
$$|(\rho-1,u,F-I)(\cdot,t)|_{H^2}^2+\int_0^t|\nabla(\rho,F)|_{H^1}^2+|\nabla u|_{H^2}^2\leq C|(\rho_0-1,u_0,F_0-I)|_{H^2}.\eqno(1.5)$$
Moreover, if $(\rho_0-1,u_0,F_0-I)\in L^1(\mathbb R^3)$, then
$$|(\rho-1,u,F-I)(t)|_{L^p}\leq C(1+t)^{-\frac{3}{2}(1-\frac{1}{p})},\ \ \forall \ \ p\in [2,6],\eqno(1.6)$$
$$|\nabla(\rho-1,u,F-I)(t)|_{H^1}\leq C(1+t)^{-\frac{5}{4}}.\eqno(1.7)$$
Finally, denote
$$(\varrho_0,m_0,\mathcal {F}_0)=(\rho_0-1,\rho_0 u_0,\rho_0 F_0-I)$$
and assume that the Fourier transform $(\hat {\varrho_0},\hat
{m}_0,\hat {\mathcal {F}_0})$ satisfies
$$|\hat{\varrho_0}|\geq c_0,\ \ \ |\hat {m_0}|\leq |\xi|^\eta,\ \ |\hat {\mathcal {F}}^T_0-\hat{\mathcal {F}}_0|\leq |\xi|^\eta,\ \ for \ \ 0\leq |\xi|\ll 1,\eqno(1.8)$$
where $c_0$ and $\eta$ are two positive constants. Then we also have
the lower bound time decay rate as
$$|(\rho-1)(t)|_{L^2}\geq c_1(1+t)^{-\frac{3}{4}},\eqno(1.9)$$
$$|u(t)|_{L^2}\geq c_1(1+t)^{-\frac{3}{4}},\eqno(1.10)$$
$$|(F-I)(t)|_{L^2}\geq c_1(1+t)^{-\frac{3}{4}},\eqno(1.11)$$
where $c_1$ is a positive constant independent of time.\\
\\{\bf Notaions.} We denote by $L^p$, $W^{m,p}$ the usual Lebesgue
and Sobolev spaces on $\mathbb R^3$ and $H^m= W^{m,2}$, with norms
$|\cdot|_{L^p}$, $|\cdot|_{W^{m,p}}$ and $|\cdot|_{H^m}$
respectively. For the sake of conciseness, we do not distinguish
functional space when scalar-valued or vector-valued functions are
involved. We denote
$\nabla=\partial_x=(\partial_1,\partial_2,\partial_3)$, where
$\partial_i=\partial_{x_i}$, $\nabla_i=\partial_i$ and put
$\partial_x^lf=\nabla^lf=\nabla(\nabla^{l-1}f)$. We assume $C$ be
a positive generic constant throughout this paper that may vary at
different places and the integration domain $\mathbb R^3$ will be
always omitted without any ambiguity. Finally, $\langle
\cdot,\cdot\rangle $ denotes the inner-product in $L^2(\mathbb
R^3)$.\\

The rest of this paper is devoted to prove Theorem 1.1. In Section
2, we first reformulate the system and  do some careful a priori
estimates for the strong solutions. Then the global existence of the
strong solutions is established by the standard continuity argument.
In Section 3 we will derive the decay-in-time estimates for the
linearized system and use the energy method to derive a
Lyapunov-type energy inequality of all the derivatives controlled by
the first order derivatives, then we utilize the decay-in-time
estimates for the linearized system to control the first order
derivatives by the higher order derivatives. Hence, the optimal
decay rates of the global strong solutions follow from these two
kinds of estimates. In section 4, we establish the lower bound time
decay rate for the global solution.

\section{Global existence}
2.1. Reformulation\\

In this subsection, we first reformulate the system (1.1). Without
loss of generality, we assume $P'(1) > 0$, and denote $\chi_0 =
(P'(1))^{-\frac{1}{2}}$ For $\rho> 0$, system (1.1) can be
rewritten as
$$
\rho_t+\rho divu+u\nabla \rho =0,\eqno(2.1a)$$$$ u^i_t+u\cdot\nabla
u^i -\frac{1}{\rho}(\mu \triangle u^i+(\lambda+\mu)\nabla_i div
u)+\frac{P'(\rho)}{\rho}\nabla_i\rho =\alpha
F^{jk}\nabla_jF^{ik},\eqno(2.1b)$$$$ F_t+u\cdot\nabla F=\nabla u
F,\eqno(2.1c)$$ where we used the condition $div(\rho F^T)=0$ for
all $t\geq 0$ which ensures that the i-th component of the vector
$div(\rho FF^T )$ is
$$\nabla_j(\rho F^{ik}F^{jk})=\rho F^{jk}\nabla_jF^{ik}+F_{ik}\nabla_j(\rho F^{jk})=\rho F^{jk}\nabla_jF^{ik}.$$

Denote $$n(t,x)=\rho(\chi_0^2t,\chi_0x)-1,\ \
v(t,x)=\chi_0u(\chi_0^2t,\chi_0x),\ \
E(t,x)=F(\chi_0^2t,\chi_0x)-I,$$ then
$$
n_t+div\  v=f-v\cdot\nabla n,\eqno(2.2a)$$$$ v^i_t -\mu \triangle
v^i-(\lambda+\mu)\nabla_i div\  v+\nabla_in-a\nabla
_jE^{ij}=g,\eqno(2.2b)$$$$ E_t-\nabla v=h-v\cdot\nabla
E,\eqno(2.2c)$$ where
$$f=-n\nabla\cdot v,\ \ h=\nabla v E,\ \ a=\frac{1}{P'(1)},$$
$$g^i=aE^{jk}\nabla_jE^{ik}-\frac{n}{1+n}(\mu
\triangle v^i+(\lambda+\mu)\nabla_i div v)-v\cdot\nabla
v^i-\left(\frac{P'(n+1)}{(1+n)P'(1)}-1\right)\nabla _in.$$ Again,
without loss of
generality, we will assume that $a=1$ for the rest of this paper.\\
\\2.2. A priori estimate\\

As a classical argument, the global existence of solutions will be
obtained by combining the local existence result with a priori
estimates. Since the local strong solutions can be proven by
standard argument of Lax-Milgram theorem and the
Schauder-Tychonoff fixed-point as [6] whose details we omit,
global solutions will follow in a standard continuity argument
after we establish (1.5) a priori. Therefore, we assume a priori
that
$$|(\rho-1,u,F-I)|_{H^2}\leq \delta_0\ll 1,\eqno(2.3)$$
which is equivalent to
$$|(n,v,E)|_{H^2}\leq \delta\ll 1.\eqno(2.4)$$
Here $\delta_0\thicksim\delta$  is small enough. This, together
with Soboles's inequality, implies in particular that
$$|(n,v,E)|_{L^\infty}\leq C\delta.\eqno(2.5)$$
This should be kept in mind in the rest of this paper.

For later use we first estimate the norm of $f,g,h$. By (2.4),
(2.5), together with Sobolev's inequality, H\"oder's inequality and
Moser-type's inequality, we easily deduce that
$$|(f,h)|_{L^2}\leq |(n,E)|_{L^\infty}|\nabla v|_{L^2}\leq C\delta|\nabla v|_{L^2},$$
$$|\nabla (f,h)|\leq |(n,E)|_{L^\infty}|\nabla^2 v|_{L^2}+|\nabla (n,E)|_{L^3}|\nabla v|_{L^6}\leq C\delta |\nabla^2v|_{L^2},$$
$$|\nabla^2(f,h)|\leq C(|(n,E)|_{L^\infty}|\nabla^3v|_{L^2}+|\nabla^2(n,E)|_{L^2}|\nabla v|_{L^\infty})\leq C\delta|\nabla^2 v|_{H^1},$$
$$\begin {array}{rl}|g|_{L^2}&\leq C(|v|_{L^3}|\nabla v|_{L^6}+|n|_{\infty}|\nabla^2v|_{L^2}+|(n,E)|_{L^\infty}|\nabla(n,E)|_{L^2})\\
&\leq C\delta(|\nabla^2v|_{L^2}+|\nabla(n,E)|_{L^2})\end
{array},$$
$$\begin {array}{rl}|\nabla g|_{L^2}\leq &C(|\nabla v|_{L^\infty}|\nabla^2v|_{L^2}+|\nabla v|_{L^\infty}|\nabla v|_{L^2}\\&+|n|_{L^\infty}|\nabla ^3v|_{L^2}
+|\nabla
n|_{L^6}|\nabla^2v|_{L^3}\\&+|(n,E)|_{L^\infty}|\nabla^2(n,E)|_{L^2}+|\nabla(n,E)|_{L^6}|\nabla(n,E)|_{L^3})\\
\leq &C\delta(|\nabla^2
v|_{H^1}+|\nabla^2(n,E)|_{L^2}),\end{array}$$ where we used the
fact $$\frac{P'(n+1)}{(1+n)P'(1)}-1\sim \mathcal {O}(1)n.$$

In what follows, a series of lemmas on the energy estimates is
given. Firstly the energy estimate of lower order for $(n,u,E)$ is
obtained in the following lemma.\\
\\{\bf Lemma 2.1.} Under the priori assumption (2.4), we have
$$\frac{1}{2}\frac{d}{dt}|(n,v,E)|^2_{L^2}+C|\nabla v|^2\leq C\delta |\nabla (n,E)|^2_{L^2}.\eqno (2.6)$$
{\bf Proof.} Multiply (2.2a), (2.2b), (2.2c) by $n,v,E$ respectively
and then integrating them over $\mathbb R^3$, we have
$$\begin {array}{rl}&\frac{1}{2}\frac{d}{dt}|(n,v,E)|^2_{L^2}+\mu|\nabla v|^2_{L^2}+(\mu+\lambda)|\nabla \cdot v|^2_{L^2}
\\=&\langle f-v\cdot\nabla n,n\rangle +\langle g,v\rangle+\langle
h-v\cdot \nabla E,E\rangle.\end{array}\eqno(2.7)$$ The three terms
on the right hand side of the above equation can be estimated as
follows.

First, it holds that
$$\begin {array}{rl}\langle f-v\cdot\nabla n,n\rangle&=\langle -n\nabla\cdot v-v\cdot\nabla n,n\rangle\\
&=\langle v\cdot\nabla n,n\rangle.\end {array}$$ It follows from
Sobolev's inequality, H\"older's inequality and (2.4) that
$$\begin {array}{rl}|\langle f-v\cdot\nabla n,n\rangle|&\leq |n|_{L^3}|v|_{L^6}|\nabla n|_{L^2}\leq C|n|_{H^1}|\nabla v|_{L^2}|\nabla n|_{L^2}
\\&\leq C\delta(|\nabla n|^2_{L^2}+|\nabla v|^2_{L^2}).\end{array}\eqno(2.8)$$
Similar to the proof of (2.8), we have $$|\langle h-v\cdot \nabla
E,E\rangle|\leq C\delta(|\nabla E|^2_{L^2}+|\nabla
v|^2_{L^2}).\eqno(2.9)$$ For the second term, we have
$$\begin {array}{rl}|\langle g^i,v^i\rangle|=&C(|\langle E^{jk}\nabla_jE^{ik},v^i\rangle|+|\langle
 \frac{n}{1+n} \triangle v^i,v^i\rangle|+|\langle
\frac{n}{1+n}\nabla div_i v,v^i\rangle|\\&+|\langle v\cdot\nabla
v^i,v^i\rangle|+|\langle(\frac{P'(n+1)}{(1+n)P'(1)}-1)\nabla
_in,v^i\rangle|).\end {array}\eqno(2.10)$$ As the proof of (2.8), it
follows from Sobolev's inequality, H\"older's inequality and (2.4)
that
$$|\langle E^{jk}\nabla_jE^{ik},v^i\rangle|\leq C\delta(|\nabla E|^2_{L^2}+|\nabla v|^2_{L^2}),\eqno(2.11)$$
$$|\langle v\cdot\nabla v,v\rangle|\leq|v|_{L^3}|v|_{L^6}|\nabla v|_{L^2}\leq C|v|_{L^6}|\nabla v|^2\leq C\delta|\nabla v|^2_{L^2},\eqno(2.12)$$
$$|\langle(\frac{P'(n+1)}{(1+n)P'(1)}-1)\nabla
_in,v^i\rangle|\leq C|n|_{L^6}|v|_{L^3}|\nabla n|_{L^2}\leq
C\delta(|\nabla n|^2_{L^2}+|\nabla v|^2_{L^2}),\eqno(2.13)$$
$$\begin {array} {rl}|\langle\frac{n}{1+n} \triangle v^i,v^i\rangle|&=|\langle\nabla (\frac{n}{1+n}) \nabla v^i,v^i\rangle|+\langle\frac{n}{1+n} \nabla v^i,\nabla v^i\rangle\\&
\leq C(|\nabla n|_{H^1}|\nabla v|_{L^2}+|n|_{L^\infty}|\nabla
v|_{L^2})\\&\leq C\delta|\nabla v|_{L^2}^2,\end {array}\eqno(2.14)$$
and similarly,
$$|\langle
\frac{n}{1+n}\nabla div v^i,v^i\rangle|\leq C\delta |\nabla
v|^2_{L^2}.\eqno(2.15)$$ Substituting (2.11)-(2.15) into (2.10)
gives that the second term is bounded by
$$|\langle g,v\rangle|\leq C\delta(|\nabla(n,E)|^2_{L^2}+|\nabla v|^2_{L^2}).\eqno(2.16)$$
Hence combining (2.7), (2.8), and (2.16) yields (2.6) since
$\delta>0$ is sufficiently small. This completes the proof of the
lemma.\\

In the following lemma we give the energy estimate of the
higher order for $(n,v,E)$.\\
\\{\bf Lemma 2.2.} Under the assumption (2.4), we have
$$\frac{1}{2}\frac{d}{dt}|\nabla(n,v,E)|^2_{H^1}+C|\nabla^2 v|^2_{H^1}\leq C\delta|\nabla(n,E)|^2_{H^1}.\eqno(2.17)$$
{\bf Proof.} Applying $\nabla$ to (2.2a), (2.2b), (2.2c) and
multiplying by $\nabla n,\nabla v,\nabla E$ respectively,
integrating over $\mathbb R^3$, we have
$$\begin {array}{rl}&\frac{1}{2}\frac{d}{dt}|\nabla(n,v,E)|^2_{L^2}+\mu|\nabla^2 v|^2_{L^2}+(\mu+\lambda)|\nabla(\nabla \cdot  v)|^2_{L^2}
\\=&\langle\nabla( f-v\cdot\nabla n),\nabla n\rangle +\langle \nabla g,\nabla v\rangle+\langle
\nabla (h-v\cdot \nabla E),\nabla E\rangle.\end{array}\eqno(2.18)$$
Now let us estimate the right-hand side term by term. First of all,
by H\"older's inequality and Sobolev's inequality, we have
$$\begin {array}{rl}|\langle \nabla f,\nabla n\rangle|+|\langle \nabla h,\nabla E\rangle|\leq &|\nabla (f,g)|_{L^2}|\nabla
(n,E)|_{L^2}\\\leq &C\delta|\nabla^2 v|_{L^2}|\nabla n|_{L^2}\\\leq&
C\delta(|\nabla^2 v|^2_{L^2}+|\nabla (n,E)|^2_{L^2}).\end {array}
$$
Next, integrating by parts, we get
$$|\langle\nabla g,\nabla v\rangle|=|\langle g,\nabla^2 v\rangle|\leq C|g|_{L^2}|\nabla^2 v|_{L^2}\leq C\delta(|\nabla
^2v|^2_{H^1}+|\nabla(n,E)|^2_{L^2}).$$ Finally by symmetry, we have
$$\begin {array}{rl}&|\langle \nabla (v\cdot\nabla n),\nabla n\rangle|+|\langle \nabla (v\cdot\nabla E),\nabla
E\rangle|\\=& |\langle \nabla v\cdot\nabla n,\nabla
n\rangle+\langle v\cdot\nabla \nabla n,\nabla n\rangle|+|\langle
\nabla v\cdot\nabla E,\nabla E\rangle+\langle  v\cdot\nabla \nabla
E,\nabla E\rangle|\\=&|\langle \nabla v\cdot\nabla n,\nabla
n\rangle+\frac{1}{2}\langle div\ v,|\nabla n|^2\rangle|+|\langle
\nabla v\cdot\nabla E,\nabla E\rangle+\frac{1}{2}\langle div\
v,|\nabla E|^2\rangle|\\ \leq &C|\nabla
v|_{L^\infty}|\nabla(n,E)|^2_{L^2}\leq C\delta(|\nabla
^2v|^2_{H^1}+|\nabla(n,E)|^2_{L^2}).\end {array}
$$
Substituting these results into (2.18), we conclude
$$\frac{1}{2}\frac{d}{dt}|\nabla(n,v,E)|^2_{L^2}+C|\nabla^2 v|^2_{L^2}\leq C\delta (|\nabla(n,E)|^2_{L^2}+|\nabla^3v|^2_{L^2}).\eqno(2.19)$$

Similarly, applying $\nabla^2$ to (2.2a), (2.2b), (2.2c) and
multiplying by $\nabla^2 n,\nabla^2 v,\nabla^2 E$ respectively,
integrating over $\mathbb R^3$, we have
$$\begin {array}{rl}&\frac{1}{2}\frac{d}{dt}|\nabla^2(n,v,E)|^2_{L^2}+\mu|\nabla^3 v|^2_{L^2}+(\mu+\lambda)|\nabla(\nabla^2 \cdot  v)|^2_{L^2}
\\=&\langle\nabla^2( f-v\cdot\nabla n),\nabla^2 n\rangle +\langle \nabla^2 g,\nabla^2 v\rangle+\langle
\nabla^2 (h-v\cdot \nabla E),\nabla^2
E\rangle.\end{array}\eqno(2.20)$$ To estimate the right-hand side of
the above equation, we note, by H\"older's inequality and Sobolev's
inequality, that
$$\begin {array}{rl}|\langle\nabla^2 f,\nabla^2 n\rangle|+|\langle\nabla^2 h,\nabla^2 E\rangle|&\leq |\nabla^2(f,h)|_{L^2}|\nabla^2(n,E)|_{L^2}
\\&\leq C\delta(|\nabla^2v|^2_{H^1}+|\nabla^2(n,E)|_{L^2}^2).
\end{array}$$
Integrating by parts, we have
$$|\langle \nabla^2g,\nabla^2v\rangle|=|\langle \nabla g,\nabla^3v\rangle|\leq |\nabla g|_{L^2}|\nabla^3 v|_{L^2}\leq C\delta(|\nabla^2v|^2_{H^1}+|\nabla^2(n,E)|^2_{L^2}).$$
Finally, by symmetry, we have
$$\begin {array}{rl}&|\langle\nabla^2(v\cdot\nabla n),\nabla^2n\rangle|\\=&|\langle\nabla^2v\cdot\nabla n,\nabla^2n\rangle
+\langle\nabla v\cdot\nabla \nabla n,\nabla^2n\rangle+\langle
v\cdot\nabla
\nabla^2n,\nabla^2n\rangle|\\=&|\langle\nabla^2v\cdot\nabla
n,\nabla^2n\rangle +\langle\nabla v\cdot\nabla \nabla
n,\nabla^2n\rangle-\frac{1}{2}\langle div\ v,|\nabla^2n|^2\rangle|\\
\leq&|\nabla^2v|_{L^6}|\nabla n|_{L^3}|\nabla^2 n|_{L^2}+|\nabla
v|_{L^\infty}|\nabla^2n|^2_{L^2}\\ \leq &C(|\nabla^3v|_{L^2}|\nabla
n|_{H^1}|\nabla^2 n|_{L^2}+|\nabla^2 v|_{H^1}|\nabla^2n|^2_{L^2})\\
\leq &C\delta(|\nabla^2v|^2_{H^1}+|\nabla^2 n|^2_{L^2}).\end
{array}$$ Similarly, we have
$$|\langle\nabla^2(v\cdot\nabla E),\nabla^2E\rangle|\leq C\delta(|\nabla^2v|^2_{H^1}+|\nabla^2 E|^2_{L^2}).$$
Putting these estimates into (2.20), we get
$$\frac{1}{2}\frac{d}{dt}|\nabla^2(n,v,E)|^2_{L^2}+C|\nabla^3 v|^2_{L^2}\leq C\delta(|\nabla^2v|^2_{H^1}+|\nabla^2 (n,E)|^2_{L^2}).\eqno(2.21)$$
Combining (2.19) and (2.21) yields (2.17) if $\delta$ is small
enough. This completes the proof of the lemma.\\

In the following lemma we give the dissipation on $|\nabla n|_{H^1}$.\\
\\{\bf Lemma 2.3.} Under the assumption (2.4), we have
$$-\frac{d}{dt}\langle div\ v,n\rangle+C|\nabla n|^2_{L^2}\leq C|\nabla v|^2_{H^1}+C\delta|\nabla E|^2_{L^2},\eqno(2.22)$$
$$\frac{d}{dt}\langle div\ v,\triangle n\rangle+C|\nabla^2 n|^2_{L^2}\leq C|\nabla v|^2_{H^2}+C\delta|\nabla^2 E|^2_{L^2}.\eqno(2.23)$$
{\bf Proof.} Notice that the condition $div(\rho F^{T})=0$ for all
$t\geq 0$ gives
$$divdiv[(1+n)(E+I)^{T}]=0,\ \ \forall\ t\geq 0.$$
Thus we have
$$\begin {array}{rl}\frac{\partial^2(
E^{ij})}{\partial x_i\partial
x_j}&=divdiv(E^T)\\&=divdiv[(1+n)(E+I)^{T}]-divdiv(nI+E^T)\\&=-\triangle
n-divdiv(nE).\end {array}\eqno(2.24)$$ Thus by applying $\nabla_i$
to (2.2b) and summing over $i$, we have
$$(div\ v)_t-(2\mu+\lambda)\triangle div\ v+2\triangle n=div\ g_1,\eqno(2.25)$$
where $$g_1=g-div(nE).$$ Multiplying the above equation by $n$, and
then integration over $\mathbb R^3$, we have
$$\begin {array}{rl}|\nabla n|^2_{L^2}=&\langle(div\ v)_t,n\rangle-
\langle(2\mu+\lambda)\triangle div\ v,n\rangle-\langle div\
g_1,n\rangle\\ =&\frac{d}{dt}\langle div\ v,n\rangle-\langle div\
v,n_t\rangle+ \langle(2\mu+\lambda)\triangle v,\nabla
n\rangle+\langle  g_1,\nabla n\rangle\\=&\frac{d}{dt}\langle div\
v,n\rangle-\langle div\ v,f\rangle+\langle div\ v,v\cdot\nabla
n\rangle+|div\ v|^2_{L^2}\\&+ \langle(2\mu+\lambda)\triangle
v,\nabla n\rangle+\langle g_1,\nabla n\rangle,\end {array}$$ where
we use the the continuity equation (2.2a). By Sobolev's, H\"older's
and Cauchy's inequalities, we obtain
$$\begin {array}{rl}&-\frac{d}{dt}\langle div\ v,n\rangle+|\nabla n|^2_{L^2}\\ \leq& C(|\nabla v|_{L^2}|f|_{L^2}
+|\nabla v|_{L^2}|v|_{L^6}|\nabla n|_{L^3}+|\nabla
v|^2_{L^2}+|\nabla^2 v|_{L^2}|\nabla n|_{L^2}\\&+|g_1|_{L^2}|\nabla
n|_{L^2})\\ \leq &C|\nabla v|^2_{H^1}+\frac{1}{2}|\nabla
n|^2_{L^2}+C\delta|\nabla(n,E)|^2_{L^2},\end {array}$$ which gives
(2.22) if $\delta$ is sufficiently small.

Multiplying (2.25) by $\triangle n$, and then integrating over
$\mathbb R^3$, we have
$$\begin {array}{rl}2|\triangle n|^2_{L^2}=&-\langle(div\
v)_t,\triangle n\rangle+ \langle(2\mu+\lambda)\triangle div\
v,\triangle n\rangle+\langle div\ g_1,\triangle n\rangle\\
=&-\frac{d}{dt}\langle div\ v,\triangle n\rangle+\langle div\
v,\triangle n_t\rangle+ \langle(2\mu+\lambda)\triangle div\
v,\triangle n\rangle\\&+\langle  div\ g_1,\triangle
n\rangle\\=&-\frac{d}{dt}\langle div\ v,\triangle n\rangle+\langle
div\ v,\triangle f\rangle-\langle div\ v,\triangle(v\cdot\nabla
n)\rangle\\&-\langle div\ v,\triangle div\ v\rangle+
\langle(2\mu+\lambda)\triangle div\ v,\triangle n\rangle+\langle
div\ g_1,\triangle n\rangle
\\=&-\frac{d}{dt}\langle div\ v,\triangle n\rangle-\langle
\nabla div\ v,\nabla f\rangle-\langle \triangle div\ v,v\cdot\nabla
n\rangle\\&+\langle \nabla div\ v,\nabla div\ v\rangle+
\langle(2\mu+\lambda)\triangle div\ v,\triangle n\rangle+\langle
div\ g_1,\triangle n\rangle.\end {array}$$ By Sobolev's, H\"older's
and Cauchy's inequalities, we have
$$\begin {array}{rl}&\frac{d}{dt}\langle div\ v,\triangle n\rangle+2|\nabla^2 n|^2_{L^2}\\ \leq& C(|\nabla^2 v|_{L^2}|\nabla f|_{L^2}
+|\nabla^3 v|_{L^2}|v|_{L^6}|\nabla n|_{L^3}+|\nabla^2
v|^2_{L^2}\\&+|\nabla^2 v|_{L^2}|\triangle n|_{L^2}+|\nabla
g_1|_{L^2}|\triangle n|_{L^2})\\ \leq &C|\nabla
v|^2_{H^2}+\frac{1}{2}|\triangle
n|^2_{L^2}+C\delta|\nabla^2(n,E)|^2_{L^2},\end {array}$$ which gives
(2.23) if $\delta$ is small enough. This completes the proof of
lemma.\\

In the following lemma we give the dissipation on $|\nabla
(E^{T}-E)|_{H^1}$.\\
\\{\bf Lemma 2.4.} Under the assumption (2.4), we have
$$-\frac{d}{dt}\langle \mathcal {W},E^{T}-E\rangle+C|\nabla (E^{T}-E)|^2_{L^2}\leq C|\nabla v|^2_{H^1}+C\delta|\nabla(n, E)|^2_{L^2},\eqno(2.26)$$
$$\frac{d}{dt}\langle \mathcal {W},\triangle (E^{T}-E)\rangle+C|\nabla^2 (E^{T}-E)|^2_{L^2}\leq C|\nabla v|^2_{H^2}+C\delta|\nabla^2 (n,E)|^2_{L^2},\eqno(2.27)$$
where $\mathcal {W}=\nabla u-(\nabla u)^T=curl \ u$.\\
\\{\bf Proof.} Taking $(2.2c)^T-(2.2c)$, we have
$$(E^{T}-E)_t+\mathcal {W}=h^T-h-v\cdot\nabla (E^T-E).\eqno(2.28)$$
Note the condition $F^{lk}\nabla_lF^{ij}=F^{lj}\nabla_lF^{ik}$ for
all $t\geq 0$, which means that
$$\nabla_kE^{ij}+E^{lk}\nabla_lE^{ij}=\nabla_jE^{ik}+E^{lj}\nabla_lE^{ik},\ \ \forall\ t\geq 0.\eqno(2.29)$$
Thus we have
$$\begin {array}{rl}&\nabla_j\nabla_kE^{ik}-\nabla_i\nabla_kE^{jk}
\\=&\nabla_k\nabla_jE^{ik}-\nabla_k\nabla_iE^{jk}\\=&\nabla_k\nabla_kE^{ij}-\nabla_k\nabla_kE^{ji}+\nabla_k(E^{lk}\nabla_lE^{ij}-E^{lj}\nabla_lE^{ik})\\
&-\nabla_k(E^{lk}\nabla_lE^{ji}-E^{li}\nabla_lE^{jk})\\=&\triangle
(E^{ij}-E^{ji})+\nabla_k(E^{lk}\nabla_lE^{ij}-E^{lj}\nabla_lE^{ik})-\nabla_k(E^{lk}\nabla_lE^{ji}-E^{li}\nabla_lE^{jk}).
\end{array}\eqno(2.30)$$
Thus by applying $curl$ to (2.2b), we have
$$\mathcal W_t-\mu\triangle\mathcal W+\triangle (E^T-E)=curl\ g+\mathcal {S},\eqno(2.31)$$
where the antisymmetric matrix $\mathcal {S}$ is defined as
$$\mathcal {S}^{ij}=\nabla_k(E^{lk}\nabla_lE^{ij}-E^{lj}\nabla_lE^{ik})-\nabla_k(E^{lk}\nabla_lE^{ji}-E^{li}\nabla_lE^{jk}).$$

Notice that the system (2.28)-(2.31) takes a similar form as the
system (2.2a)-(2.25). Thus after a similar argument as Lemma 2.3,
(2.26) and (2.27) follows. The proof of lemma is completed.\\

Finally, in the following lemma we give the dissipation on $|\nabla
E|_{H^1}$.\\
\\{\bf Lemma 2.5.} Under assumption (2.4), we have
$$|\nabla E|^2_{L^2}\leq C|\nabla(n,E^T-E)|_L^2,\eqno(2.32)$$
$$|\nabla^2 E|^2_{L^2}\leq C|\nabla^2(n,E^T-E)|_L^2.\eqno(2.33)$$
\\{\bf Proof.} Combining (2.24) and (2.30), we have
$$\begin {array}{rl}\triangle div\ E&=\nabla div div\ E-curl curl div\ E
\\&=-\triangle\nabla n-\nabla div div(nE)+\triangle curl(E-E^T)+curl\mathcal {S}.\end {array}\eqno(2.34)$$
Thus using the property of Riesz potential, (2.4) and (2.5), we
arrive at
$$\begin {array}{rl}|div E|^2_{L^2}&\leq C(|\nabla n|^2_{L^2}+|\nabla (E^T-E)|^2_{L^2}+|\nabla(nE)|^2_{L^2}+|E\nabla E|^2_{L^2})
\\&\leq C|\nabla (n,E^T-E)|^2_{L^2}+C\delta|\nabla^2 E|^2_{L^2},\end {array}$$
and
$$\begin {array}{rl}|\nabla div E|^2_{L^2}&\leq C(|\nabla^2 n|^2_{L^2}+|\nabla^2 (E^T-E)|^2_{L^2}+|\nabla^2(nE)|^2_{L^2}+|\nabla (E\nabla E)|^2_{L^2})
\\&\leq C|\nabla^2 (n,E^T-E)|^2_{L^2}+C\delta|\nabla E|^2_{L^2}.\end {array}$$
Under the above estimate, we may deduce from (2.29) that
$$\begin {array}{rl}|\nabla E|^2_{L^2}&\leq |div E|^2_{L^2}+|curl E|^2_{L^2}
\\&\leq C|\nabla (n,E^T-E)|^2_{L^2}+C\delta|\nabla E|^2_{L^2}+|E\nabla E|^2_{L^2}\\ &\leq C|\nabla (n,E^T-E)|^2_{L^2}+C\delta|\nabla E|^2_{L^2},\end {array}$$
and
$$\begin {array}{rl}|\nabla^2 E|^2_{L^2}&\leq |\nabla div E|^2_{L^2}+|\nabla curl E|^2_{L^2}
\\&\leq C|\nabla^2 (n,E^T-E)|^2_{L^2}+C\delta|\nabla^2 E|^2_{L^2}+|\nabla (E\nabla E)|^2_{L^2}\\ &\leq C|\nabla^2 (n,E^T-E)|^2_{L^2}+C\delta|\nabla^2 E|^2_{L^2}.\end {array}$$
This proves (2.32) and (2.33), and the proof lemma is completed.\\

Now we are in a position to verify (2.4). Since $\delta>0$ is
sufficiently small, from Lemma 2.1-Lemma 2.5, we can choose a
constant $D_1>0$ suitably large such that
$$\begin {array}{rl}&\frac{d}{dt}\{D_1|(n,v,E)|^2_{H^2}+\langle div\ v,\triangle n-n\rangle
+\langle \mathcal {W},\triangle (E^T-E)-(E^T-E)\rangle\}\\
&+C(|\nabla (n,E)|^2_{H^1}+|\nabla v|^2_{H^2})\leq 0,\end
{array}$$ for any $t\geq 0$, which implies
$$|(n,v,E)|^2_{H^2}+\int_0^t(|\nabla (n,E)|^2_{H^1}+|\nabla v|^2_{H^2})\leq C|(n_0,v_0,E_0)|^2_{H^2},\eqno(2.35)$$
since
$$D_1|(n,v,E)|^2_{H^2}+\langle div\ v,\triangle n-n\rangle
+\langle \mathcal {W},\triangle (E^T-E)-(E^T-E)\rangle\sim
|(n,v,E)|^2_{H^2}.$$ Then (2.35) gives (2.4). Thus we prove the
global existence result of Theorem 1.1.

\section{Convergence rate of the solution}

$\ \ \ \ $In this section we shall prove the decay rates of the
solution to finish the proof of Theorem 1.1. In Section 3.1, we list
some elementary conclusion on the decay-in-time estimates for the
linearized system and a useful inequality. In Section 3.2, we shall
first obtain the energy inequality for the derivatives of the orders
from the first to the third, and then we show a decay-in-time
estimate for the first order derivatives, where the error is related
to the derivatives of the higher order. Finally, by combining these
estimates we get the optimal decay rates.\\
\\3.1. Spectral analysis and linear $L^2$ estimates\\

We first note that the linearized system (2.2a)-(2.25) depends only
on $(n,div\ v)$ while the linearized system (2.28)-(2.31) also
depends only on $(\mathcal {W},E^T-E)$. Denote by $\Lambda^s$ the
pseudo differential operator defined by
$$\Lambda^su=\mathscr{F}^{-1}(|\xi|^s\hat u(\xi)),$$ and let
$$m=\Lambda^{-1}div\ v$$
be the ``compressible part" of the velocity, and
$$\omega=\Lambda^{-1}\mathcal {W}=\Lambda^{-1}curl\ v$$
be the ``incompressible part" of the velocity. We finally obtain
$$
n_t+\Lambda d=f-v\cdot\nabla n,\eqno(3.1a)$$
$$d_t-(2\mu+\lambda)\triangle d-2\Lambda n=\Lambda^{-1}div\ g_1,\eqno(3.1b)$$
$$(E^{T}-E)_t+\Lambda \omega=h^T-h-v\cdot\nabla (E^T-E),\eqno(3.2a)$$
$$\omega_t-\mu\triangle\omega-\Lambda (E^T-E)=\Lambda^{-1}curl\ g+\Lambda^{-1}\mathcal {S}.\eqno(3.2b)$$
Indeed, as the definition of $d$ and $\omega$, and the relation
$$v=-\Lambda^{-1}\nabla d+\Lambda^{-1}curl\  \omega$$
involve pseudo-differential operators of degree zero, the estimates
in space $H^{l}(\mathbb{R}^{3})$ for the original function $v$ will
be the same as for $(d,\omega)$.

Here, we just discuss the system (3.1) for example, since the system
(3.2) is the same as system (3.1). To use the $L^p-L^q$ estimates of
the linear problem for the nonlinear system (3.1) and system (3.2),
we rewrite the solution of (3,1) as
$$U(t)=K(t)U_0+\int_0^tK(t-\tau)G(\tau)d\tau\ \ \ t\geq0,$$
where we use the notations
$$U=[n,d]^T,\ \ U_0=[n_0,d_0]^T,\ \ G=[f-v\cdot\nabla n,\Lambda^{-1}div\ g_1]^T,$$
and $K(t)$ is the solution semigroup defined by $K(t)=e^{tB},\ t\geq
0$, with $B$ being a matrix-valued differential operator given by
$$
B:=\left(\begin{array}{cc}0&-\Lambda\\2\Lambda&(2\mu+\lambda)\triangle
\end{array}\right) .$$

Now we aim to analyze the differential operator $B$ in terms of
its Fourier expression $A$ and to show the long time properties of
the semigroup $K(t)$. For this purpose, we need to consider the
following linearized system
$$U_t=BU.\eqno(3.3)$$
Applying the Fourier transform to system (3.3), we have
$$\partial_t\hat U=A(\xi)\hat U,\ \ \hat U(0)=\hat U_0,$$
where $\hat U(t)=\hat
U(\xi,t)=\mathscr{F}U(\xi,t),\xi=(\xi_1,\xi_2,\xi_3)$ and $A(\xi)$
is defined as
$$
A(\xi):=\hat
B=\left(\begin{array}{cc}0&-|\xi|\\2|\xi|&-(2\mu+\lambda)|\xi|^2
\end{array}\right) .$$
The characteristic polynomial of $A(\xi)$ is
$\kappa^2+(2\mu+\lambda)\kappa+2|\xi|^2$, which implies the
eigenvalues are
$$\kappa_{\pm }=-(\mu+\frac{1}{2}\lambda)|\xi|^2\pm\frac{1}{2}i\sqrt{8|\xi|^2-(2\mu+\lambda)^2|\xi|^4}.$$
The semigroup $e^{tA}$ is expressed as
$$\begin {array}{rl}e^{tA}&=e^{\kappa_+t}\frac{A(\xi)-\kappa_-I}{\kappa_+-\kappa_-}+e^{\kappa_-t}\frac{A(\xi)-\kappa_+I}{\kappa_--\kappa_+}
\\&=\left(\begin{array}{cc}\frac{\kappa_+e^{\kappa_-t}-\kappa_-e^{\kappa_+t}}{\kappa_+-\kappa_-}&-\frac{|\xi|(e^{\kappa_+t}-e^{\kappa_-t})}{\kappa_+-\kappa_-}
\\\frac{2|\xi|(e^{\kappa_+t}-e^{\kappa_-t})}{\kappa_+-\kappa_-}&\frac{\kappa_+e^{\kappa_-t}-\kappa_-e^{\kappa_+t}}{\kappa_+-\kappa_-}-\frac{(2\mu+\lambda)|\xi|^2(e^{\kappa_+t}-e^{\kappa_-t})}{\kappa_+-\kappa_-}\end{array}\right)
.\end{array}$$ Thus the semigroup $K(t)$ has the following
properties on the decay in time, which can be found in [10,11,14].\\
\\{\bf Lemma 3.1.} Let $k\geq 0$ be an integer and $1\leq l\leq
2$.  Then for any $t\geq 0$, the solution $U(t)=(n(t),d(t))$ of
system (3.6) satisfies
$$|\nabla^kK(t)U_0|_{L^2}\leq C(1+t)^{-\sigma(l,2;k)}(|\hat U_0|_{L^\frac{l}{l-1}}+|\nabla^kU_0|_{L^2})\leq C(1+t)^{-\sigma(l,2;k)}|(n,v)|_{L^l\cap H^k},$$
where the decay rate is measured by
$$\sigma(l,2;k)=\frac{3}{2}(\frac{1}{l}-\frac{1}{2})+\frac{k}{2}.\eqno(3.4)$$
Moreover, if $|\hat n_0|\geq c_0,\hat d_0=0$ for $0\leq|\xi|\ll 1$,
then there exists a positive constant $c_2$ such that
$$|n(t)|_{L^2}\geq c_2(1+t)^{-\frac{3}{4}},$$
$$|d(t)|_{L^2}\geq c_2(1+t)^{-\frac{3}{4}}.$$
Finally, if $|(\hat n_0,\hat d_0)|\leq |\xi|^\eta$ for
$0\leq|\xi|\ll 1$, then there exists a positive constant $C$ such
that
$$|(n,d)(t)|_{L^2}\leq {C(1+t)}^{-\frac{\eta}{2}-\frac{3}{4}}|(n_0,d_0)|_{L^2}.$$\\

We finish this subsection by listing an elementary but useful
inequality [4]:\\
\\{\bf Lemma 3.2.} If $r_1>1$ and $r_2\in[0,r_1]$, then it holds that
$$\int_0^\tau(1+t-\tau)^{-r_1}(1+\tau)^{-r_2}\leq C(r_1,r_2)(1+t)^{-r_2}.$$\\
3.2. Convergence rates\\

Now we will show the energy inequality as follows:\\
\\{\bf Lemma 3.3.} Under the assumption (2.4), let $(n,v,E)$ be the
solution to the initial value problem (2.2), then there are two
positive constants $C$ and $D_2$ such that if $\delta>0$ in (2.4) is
small enough, it holds
$$\frac{d}{dt}M(t)+D_2M(t)\leq C|\nabla (n,v,E)(t)|^2_{L^2},\eqno(3.5)$$
where the energy function $M(t)$ defined by (3.7) is equivalent to
$|\nabla (n,v,E)|^2_{H^1}$, that is, there exists a positive
constant $C_1>0$ such that $$\frac{1}{C_1}|\nabla
(n,v,E)(t)|^2_{H^1}\leq M(t)\leq C_1|\nabla (n,v,E)(t)|^2_{H^1}.$$
{\bf Proof.} Since $\delta>0$ is sufficiently small, from Lemma
2.2-Lemma 2.5, we can choose a constant $D_2>0$ suitably large such
that
$$\begin {array}{rl}&\frac{d}{dt}\{D_2|\nabla(n,v,E)|^2_{H^1}+\langle div\ v,\triangle n\rangle
+\langle \mathcal {W},\triangle (E^T-E)\rangle\}\\
&+C(|\nabla^2 (n,E)|^2_{L^2}+|\nabla^2 v|^2_{H^1})\leq C\delta
|\nabla(n,v,E)|_{L^2}.\end {array}\eqno(3.6)$$ Define the energy
functional
$$M(t)=D_2|\nabla(n,v,E)|^2_{H^1}+\langle div\ v,\triangle n\rangle
+\langle \mathcal {W},\triangle (E^T-E)\rangle,\eqno(3.7)$$ for any
$t\geq 0$, where it is noticed that $M(t)$ is equivalent to
$|\nabla(n,v,E)|^2_{H^1}$ since $D_2$ can be large enough. Adding
$|\nabla (n,v,E)|$ to both sides of (3.6) gives (3.5). This
completes the proof of the lemma.\\

If we define
$$N(t)=\sup_{0\leq\tau\leq t}(1+\tau)^{\frac{5}{2}}M(\tau),\eqno(3.8)$$
then
$$|\nabla (n,v,E)(t)|_{H^1}\leq C\sqrt{M(t)}\leq C(1+t)^{-\frac{5}{4}}\sqrt{N(t)}.\eqno(3.9)$$
To close the estimate (3.5), we shall estimate the decay rate of the
first order derivatives, this will be based on Lemma 3.1 about the
decay estimates
on the semigroup $K(t)$. Precisely, we have the following lemma.\\
\\{\bf Lemma 3.4.} Under the assumption (2.4), let $(n,v,E)$ be
the solution to the initial value problem (2.2). Then we have
$$|\nabla (n,v,E)(t)|_{L^2}\leq C(1+t)^{-\frac{5}{4}}(K_0+\delta\sqrt{N(t)}),\eqno(3.10)$$
where $K_0=|(n_0,v_0,E_0)|_{L^1\cap H^2}$.\\
\\{\bf Proof.} From the Duhamel's principle, it holds that
$$
\left(\begin{array}{c}n\\d
\end{array}\right) =K(t)\left(\begin{array}{c}n_0\\d_0
\end{array}\right)+\int K(t-\tau)G(\tau) d\tau.$$
Thus from Lemma 3.1, we have
$$|(n,d)|\leq CK_0(1+t)^{-\sigma(1,2;1)}+C\int_0^t(1+t-\tau)^{-\sigma(1,2;1)}(|\hat G(\tau)|_{L^\infty}+|G(\tau)|_{H^1})d\tau,\eqno(3.11)$$
where $\sigma(1,2;1)=\frac{5}{4}$ by (3.4). By H\"older's inequality
and Sobolev's inequality, the nonlinear source terms can be
estimated as follows:
$$|\hat G(t)|_{L^\infty}\leq C|(f-v\cdot \nabla n,g_1)|_{L^1}\leq C\delta(|\nabla(n,E)|_{L^2}+|\nabla v|_{H^1}),$$
$$|G(t)|_{H^1}\leq C|(f-v\cdot \nabla n,g_1)|_{H^1}\leq C\delta|\nabla(n,v,E)|_{H^1}+C|\nabla n|_{H^1}|\nabla^3 v|_{L^2}).$$
Putting these estimates into (3.11), by (2.35), (3.9), Lemma 3.2 and
H\"older's inequality, we arrive at
$$\begin {array}{rl}|\nabla(n,d)(t)|_{L^2}\leq &CK_0(1+t)^{-\frac{5}{4}}+C\delta\int_0^t(1+t-\tau)^{-\frac{5}{4}}
(1+\tau)^{-\frac{5}{4}}\sqrt{N(\tau)}d\tau\\&+C\int_0^t(1+t-\tau)^{-\frac{5}{4}}
(1+\tau)^{-\frac{5}{4}}\sqrt{N(\tau)}|\nabla^3v(\tau)|_{L^2}d\tau\\
\leq
&C(1+t)^{-\frac{5}{4}}(K_0+\delta\sqrt{N(t)})\\&+C\sqrt{N(t)}\int_0^t(1+t-\tau)^{-\frac{5}{4}}
(1+\tau)^{-\frac{5}{4}}|\nabla^3v(\tau)|_{L^2}d\tau\\ \leq
&C(1+t)^{-\frac{5}{4}}(K_0+\delta\sqrt{N(t)})\\&+C\sqrt{N(t)}(\int_0^t(1+t-\tau)^{-\frac{5}{2}}
(1+\tau)^{-\frac{5}{2}}d\tau)^{\frac{1}{2}}(\int_0^t|\nabla^3v(\tau)|^2_{L^2}d\tau)^{\frac{1}{2}}
\\ \leq
&C(1+t)^{-\frac{5}{4}}(K_0+\delta\sqrt{N(t)}).
\end {array}$$
Similarly, we have $$|\nabla(\omega, E^T-E)(t)|_{L^2}\leq
C(1+t)^{-\frac{5}{4}}(K_0+\delta\sqrt{N(t)}).$$ Combining the above
two inequalities, Lemma 2.5, and the relation of $v$ and
$(d,\omega)$, we get (3.10). This completes the proof of the
lemma.\\

Now we are in a position to prove (1.6)-(1.7) in Theorem 1.1.
Applying the Gronwall's inequality to the Lyapunov-type inequality
(3.5), by (3.10), we get
$$\begin {array}{rl}M(t)&\leq M(0)e^{-D_2t}+C\int_0^te^{-D_2(t-\tau)}|\nabla(n,v,E)|^2_{L^2}d\tau
\\&\leq M(0)e^{-D_2t}+C\int_0^te^{-D_2(t-\tau)}(1+\tau)^{-\frac{5}{2}}(K^2_0+\delta^2N(\tau))d\tau
\\&\leq C(1+t)^{-\frac{5}{2}}(M(0)+K^2_0+\delta^2N(t)).\end {array}$$
In view of (3.8), we have
$$N(t)\leq C(M(0)+K^2_0+\delta^2N(t)),$$
which implies
$$N(t)\leq C(M(0)+K_0^2)\leq CK_0^2,$$
since $\delta>0$ is sufficiently small. Thus (3.9) gives
$$|\nabla(n,v,E)(t)|_{H^1}\leq CK_0(1+t)^{-\frac{5}{4}}.\eqno (3.12)$$
This proves (1.7). Now for (1.6), first by Sovolev's inequality and
(3.12), we have
$$|(n,v,E)(t)|_{L^6}\leq C|\nabla(n,v,E)(t)|_{L^2}\leq CK_0(1+t)^{-\frac{5}{4}}.\eqno(3.13)$$
Meanwhile, using Lemma 3.1 and Lemma 3.2, it follows from the
Duhamel's principle that $$\begin {array}{rl}&|(n,d)(t)|_{L^2}\\
\leq &
CK_0(1+t)^{-\frac{3}{4}}+C\int_0^t(1+t-\tau)^{-\frac{3}{4}}|(f-v\cdot\nabla
n,g_1)(\tau)|_{L^1\cap L^2}d\tau\\ \leq&
CK_0(1+t)^{-\frac{3}{4}}+C\delta\int_0^t(1+t-\tau)^{-\frac{3}{4}}|\nabla(n,v,E)(\tau)|_{H^1}d\tau
\\ \leq&
CK_0(1+t)^{-\frac{3}{4}}+C\delta
K_0\int_0^t(1+t-\tau)^{-\frac{3}{4}}(1+\tau)^{-\frac{5}{4}}d\tau
\\ \leq&
CK_0(1+t)^{-\frac{3}{4}}.
\end {array}\eqno(3.14)$$
Similarly, we have
$$|(\omega,E^T-E)(t)|_{L^2}\leq CK_0(1+t)^{-\frac{3}{4}}.\eqno(3.15)$$
Finally, we derive the time-decay-rate on $|E(t)|_{L^2}$. From
(2.29) and (2.34), we have
$$|\Lambda^{-1}curl\ E(t)|_{L^2}\leq C|\Lambda^{-1}(E\nabla E)(t)|_{L^2}\leq C|E(t)\nabla E(t)|_{L^{\frac{6}{5}}}\leq CK_0(1+t)^{-\frac{5}{4}},\eqno(3.16)$$
$$\begin {array}{rl}|\Lambda^{-1} div\ E(t)|_{L^2}&\leq C(|(n,E-E^T)(t)|_{L^2}+|n(t)E(t)|_{L^2}+|\Lambda^{-2}\mathcal {S}(t)|_{L^2})\\&
\leq CK_0(1+t)^{-\frac{3}{4}}+C|\Lambda^{-1}(E\nabla
E)(t)|_{L^2}\\&\leq CK_0(1+t)^{-\frac{3}{4}}.\end {array}$$ The
above two inequalities give
$$|E(t)|_{L^2}\leq CK_0(1+t)^{-\frac{3}{4}}.\eqno(3.17)$$
Combining (3.14), (3.15), (3.17) and the relation of $v$ and
$(d,\omega)$, we have
$$|(n,v,E)(t)|_{L^2}\leq CK_0(1+t)^{-\frac{3}{4}}.\eqno(3.18)$$
Hence, by the interpolation, it follows from (3.13), (3.18) that for
any $2\leq p\leq 6$
$$|(n,v,E)(t)|_{L^p}\leq |(n,v,E)(t)|^\theta_{L^2}|(n,v,E)(t)|^{1-\theta}_{L^6} \leq C(1+t)^{-\frac{3}{2}(1-\frac{1}{p})},$$
where $\theta=\frac{6-p}{2p}$. The proof of (1.6)-(1.7) is
completed.

\section{Lower bound time decay rate}

$\ \ \ \ $In this section, we investigate the lower bound time
decay for global solutions. Define
$$\varrho(t,x)=\rho(t,x)-1,\ \ m(t,x)=\rho u,\ \ \mathcal F=\rho F-I.$$
Then the condition $div\mathcal F^T=0$ ensures that
$$\begin {array}{rl}div(\rho FF^T)&=div[(\mathcal F+I)(\frac{\mathcal F^T+I}{\varrho+1})]\\&
=div\mathcal F-\nabla \varrho+div(\frac{-\varrho\mathcal F+\mathcal
F\mathcal F^T-\varrho\mathcal F^T+\varrho^2I}{1+\varrho}),\end
{array}$$ and
$$divdiv\mathcal F=divdiv\mathcal F^T=0.$$ Thus we have
the following system which only depends on $(\varrho,div\ m)$
$$\varrho_t+div\ m=0,$$
$$(div\ m)_t-(2\mu+\lambda)\triangle (div\ m)+(1+\alpha)\triangle\varrho=G_1,$$
where $$\begin {array}{rl}G_1=&\alpha divdiv(\frac{-\varrho\mathcal
F+\mathcal F\mathcal F^T-\varrho\mathcal
F^T+\varrho^2I}{1+\varrho})+\triangle(\varrho-P(1+\varrho))\\&+(2\mu+\lambda)\triangle
(\frac{\varrho div\ m}{1+\varrho})-divdiv(\rho u\otimes
u).\end{array}$$ By H\"older's inequality and Sobolev's inequality,
it is easy to verify that
$$(|\widehat{\Lambda^{-2} G_1(t)}|_{L^\infty}+|\Lambda^{-1} G_1(t)|_{L^2})\leq C|(\varrho,m,\mathcal F)|^2_{H^2}\leq C(1+t)^{-\frac{3}{2}}.$$
Thus by Duhamel's principle, Lemma 3.1 and the condition (1.8), we
have
$$\begin {array}{rl}&|(\varrho,\Lambda^{-1}div\ m)(t)|_{L^2}\\ \geq &|K(t)(\varrho_0,\Lambda^{-1}div\ m_0)|_{L^2}-\int_0^t|K(t-\tau)(0,\Lambda^{-1}G_1(\tau))|_{L^2}d\tau
\\ \geq&
c_1(1+t)^{-\frac{3}{4}}-C\int_0^t(1+t-\tau)^{-\frac{5}{4}}(|\Lambda^{-2}\hat
G_1(\tau)|_{L^\infty}+|\Lambda^{-1} G_1(\tau)|_{L^2})d\tau\\ \geq&
c_1(1+t)^{-\frac{3}{4}}-C\int_0^t(1+t-\tau)^{-\frac{5}{4}}(1+\tau)^{-\frac{3}{2}}d\tau\\
\geq& c_2(1+t)^{-\frac{3}{4}}.
\end
{array}\eqno(4.1)$$ Hence (1.9) is proved.

On the other hand, the condition
$F^{lk}\nabla_lF^{ij}=F^{lj}\nabla_lF^{ik}$ means that
$$\begin{array}{rl}&\nabla_k\mathcal F^{ij}+\mathcal F^{lk}\nabla_l(\frac{\mathcal F^{ij}+\delta_{ij}}{1+\varrho})
+\nabla_k(\frac{-\varrho\mathcal
F^{ij}+\delta_{ij}}{1+\varrho})\\=&\nabla_j\mathcal F^{ik}+\mathcal
F^{lj}\nabla_l(\frac{\mathcal F^{ik}+\delta_{ik}}{1+\varrho})
+\nabla_j(\frac{-\varrho\mathcal
F^{ik}+\delta_{ik}}{1+\varrho}),\end{array}$$ where
$$\delta=
\left\{
\begin{array}{l}
0\ if\ i\neq j,\\
1\ if\ i=j.
\end{array}
\right.
$$
Using the fact $div\ \mathcal F^T=0$, we have
$$\begin {array}{rl}&\nabla_j\nabla_k\mathcal F^{ik}-\nabla_i\nabla_k\mathcal F^{jk}
\\=&\nabla_k\nabla_j\mathcal F^{ik}-\nabla_k\nabla_i\mathcal F^{jk}\\=&\nabla_k\nabla_k\mathcal F^{ij}-\nabla_k\nabla_k\mathcal F^{ji}
+\nabla_k(\mathcal F^{lk}\nabla_l(\frac{\mathcal
F^{ij}+\delta_{ij}}{1+\varrho}) +\nabla_k(\frac{-\varrho\mathcal
F^{ij}+\delta_{ij}}{1+\varrho}))\\&-\nabla_k(\mathcal
F^{lj}\nabla_l(\frac{\mathcal F^{ik}+\delta_{ik}}{1+\varrho})
+\nabla_j(\frac{-\varrho\mathcal F^{ik}+\delta_{ik}}{1+\varrho}))
\\&-\nabla_k(\mathcal F^{lk}\nabla_l(\frac{\mathcal
F^{ji}+\delta_{ji}}{1+\varrho}) +\nabla_k(\frac{-\varrho\mathcal
F^{ji}+\delta_{ji}}{1+\varrho})) \\&+\nabla_k(\mathcal
F^{li}\nabla_l(\frac{\mathcal F^{jk}+\delta_{jk}}{1+\varrho})
+\nabla_k(\frac{-\varrho\mathcal F^{jk}+\delta_{jk}}{1+\varrho}))
\\=&\nabla_k\nabla_k\mathcal F^{ij}-\nabla_k\nabla_k\mathcal F^{ji}
+\nabla_k(\mathcal F^{lk}\nabla_l(\frac{\mathcal F^{ij}}{1+\varrho})
+\nabla_k(\frac{-\varrho\mathcal
F^{ij}}{1+\varrho}))\\&-\nabla_k(\mathcal
F^{lj}\nabla_l(\frac{\mathcal F^{ik}}{1+\varrho})
+\nabla_j(\frac{-\varrho\mathcal F^{ik}}{1+\varrho}))
\\&-\nabla_k(\mathcal F^{lk}\nabla_l(\frac{\mathcal
F^{ji}}{1+\varrho}) +\nabla_k(\frac{-\varrho\mathcal
F^{ji}}{1+\varrho})) \\&+\nabla_k(\mathcal
F^{li}\nabla_l(\frac{\mathcal F^{jk}}{1+\varrho})
+\nabla_k(\frac{-\varrho\mathcal
F^{jk}}{1+\varrho}))\\&-\nabla_i(\mathcal
F^{lj}\nabla_l\frac{1}{1+\varrho})+\nabla_j(\mathcal
F^{li}\nabla_l\frac{1}{1+\varrho})
\\=&\triangle(\mathcal F^{ij}-\mathcal F^{ji})
+\nabla_k(\nabla_l(\mathcal F^{lk}\frac{\mathcal F^{ij}}{1+\varrho})
+\nabla_k(\frac{-\varrho\mathcal
F^{ij}}{1+\varrho}))\\&-\nabla_k(\nabla_l(\mathcal
F^{lj}\frac{\mathcal F^{ik}}{1+\varrho})
+\nabla_j(\frac{-\varrho\mathcal F^{ik}}{1+\varrho}))
\\&-\nabla_k(\nabla_l(\mathcal F^{lk}\frac{\mathcal
F^{ji}}{1+\varrho}) +\nabla_k(\frac{-\varrho\mathcal
F^{ji}}{1+\varrho})) \\&+\nabla_l(\nabla_k(\mathcal
F^{li}\frac{\mathcal F^{jk}}{1+\varrho})
+\nabla_k(\frac{-\varrho\mathcal
F^{jk}}{1+\varrho}))\\&+\nabla_i(\mathcal \nabla_l\frac{\varrho
F^{lj}}{1+\varrho})-\nabla_j(\nabla_l\frac{\varrho\mathcal
F^{li}}{1+\varrho}).
\end{array}$$
Thus by applying $curl$ to (1.1b) we have
$$(curl\ m)_t-\mu\triangle(curl\ m)+\alpha\triangle (\mathcal F^T-\mathcal F)=H_1\eqno(4.2)$$
where
$$\begin{array}{rl}H_1^{ij}=&\nabla_j(\mu\triangle\frac{\varrho m}{1+\varrho})-\nabla_i(\mu\triangle\frac{\varrho m}{1+\varrho})+\nabla_j(div(\rho u\otimes
u))-\nabla_i(div(\rho u\otimes u))\\&+\nabla_k(\nabla_l(\mathcal
F^{lk}\frac{\mathcal F^{ij}}{1+\varrho})
+\nabla_k(\frac{-\varrho\mathcal
F^{ij}}{1+\varrho}))-\nabla_k(\nabla_l(\mathcal F^{lj}\frac{\mathcal
F^{ik}}{1+\varrho}) +\nabla_j(\frac{-\varrho\mathcal
F^{ik}}{1+\varrho}))
\\&-\nabla_k(\nabla_l(\mathcal F^{lk}\frac{\mathcal
F^{ji}}{1+\varrho}) +\nabla_k(\frac{-\varrho\mathcal
F^{ji}}{1+\varrho})) +\nabla_l(\nabla_k(\mathcal
F^{li}\frac{\mathcal F^{jk}}{1+\varrho})
+\nabla_k(\frac{-\varrho\mathcal
F^{jk}}{1+\varrho}))\\&+\nabla_i(\mathcal \nabla_l\frac{\varrho
F^{lj}}{1+\varrho})-\nabla_j(\nabla_l\frac{\varrho\mathcal
F^{li}}{1+\varrho}).\end{array}$$ We also note that
$$\begin {array}{rl}\nabla\times(u\times\rho F^T)&=\rho F^T\cdot\nabla u-\rho F^Tdiv\ u-u\cdot\nabla(\rho F^T)+udiv(\rho F^T)
\\&=\rho F^T(\nabla u)^T-\rho F^Tdiv\ u-(u\cdot\nabla\rho)F^T-\rho (u\cdot\nabla F^T)
\\&=\rho (F^T(\nabla u)^T-u\cdot\nabla F^T)-(\rho div\ u+u\cdot\nabla\rho)F^T\\&
=\rho F^T_t+\rho_t  F^T=(\rho F^T)_t,\end {array}$$ where we used
the condition $div(\rho F^T)=0$. Thus we have
$$(\mathcal F^T-\mathcal F)_t+curl\ m=H_2,\eqno(4.3)$$
where $$H_2=\nabla\times(u\times\mathcal F^T)+curl\ (\varrho
u)-(\nabla\times(u\times\mathcal F^T)+curl\ (\varrho u))^T.$$ By
H\"older's inequality and Sobolev's inequality, it is easy to verify
that
$$(|(\widehat{\Lambda^{-2}H_1},\widehat{\Lambda^{-1} H_2})(t)|_{L^\infty}+| (\Lambda^{-1}H_1,H_2)(t)|_{L^2})\leq C|(\varrho,m,\mathcal F)|^2_{H^2}\leq C(1+t)^{-\frac{3}{2}}.$$
Duhamel's principle, Lemma 3.1 and the condition (1.8), we have the
following estimates of system (4.2)-(4.3):
$$\begin {array}{rl}&|(\mathcal F^T-\mathcal F,\Lambda^{-1}curl\ m)(t)|_{L^2}\\ \leq &|K(t)(\mathcal F_0^T-\mathcal F_0,\Lambda^{-1}curl\ m_0)|_{L^2}-\int_0^t|K(t-\tau)(H_2(\tau),\Lambda^{-1}H_1(\tau)|_{L^2}d\tau
\\ \leq&
C+C\int_0^t(1+t-\tau)^{-\frac{5}{4}}\\&\times(|(\widehat{\Lambda^{-2}H_1},\widehat{\Lambda^{-1}
H_2})(\tau)|_{L^\infty}+| (\Lambda^{-1}H_1,H_2)(\tau)|_{L^2})d\tau\\
\leq&
C(1+t)^{-\frac{\eta}{2}-\frac{3}{4}}+C\int_0^t(1+t-\tau)^{-\frac{5}{4}}(1+\tau)^{-\frac{3}{2}}d\tau\\
\leq& C(1+t)^{-\min(\frac{\eta}{2}+\frac{3}{4},\frac{5}{4})}.
\end
{array}\eqno(4.4)$$ Combining (4.1) and (4.4) gives
$$\begin {array}{rl}|m(t)|_{L^2}&\geq |\Lambda^{-1}div\ m(t)|_{L^2}-|\Lambda^{-1}curl\ m(t)|_{L^2}
\\&\geq c_3(1+t)^{-\frac{3}{4}}-C(1+t)^{-\min(\frac{\eta}{2}+\frac{3}{4},\frac{5}{4})}\\&\geq c_4(1+t)^{-\frac{3}{4}}.\end {array}$$
Hence (1.10) is proved.

By (4.4), we also have
$$|(E^T-E)(t)|_{L^2}\leq C|(\mathcal F^T-\mathcal F)(t)|_{L^2}\leq C(1+t)^{-\min(\frac{\eta}{2}+\frac{3}{4},\frac{5}{4})}.\eqno (4.5)$$
Thus from (2.34), (4.1) and (4.5) , we obtain
$$\begin {array}{rl}|\Lambda^{-1}div\ E(t)|_{L^2}&\geq |n(t)|_{L^2}-|n(t)E(t)|_{L^2}-|(E^T-E)(t)|_{L^2}-|\Lambda^{-2}\mathcal {S}(t)|_{L^2}
\\&\geq c_2(1+t)^{-\frac{3}{4}}-C(1+t)^{-\min(\frac{\eta}{2}+\frac{3}{4},\frac{5}{4})}\\&\geq c_5(1+t)^{-\frac{3}{4}}.\end {array}$$
Combining the above inequality with (3.16), we arrive at
$$|E(t)|_{L^2}\geq |\Lambda^{-1}div\ E(t)|_{L^2}-|\Lambda^{-1}curl\ E(t)|_{L^2}\geq c_6(1+t)^{-\frac{3}{4}}.$$
Thus, (1.11) is proved and this completes the proof of Theorem
1.1.\\
\\{\bf Acknowledgments}\\

Xianpeng Hu's research was supported in part by the National
Science Foundation. Guochun Wu's research was supported by China
Scholarship Council (File No. 201206310033).

\end{document}